\documentclass[10pt,a4paper]{amsart}
      \usepackage[english]{babel}
      \usepackage{amsmath}
      \usepackage{amstext}
      \usepackage{amssymb}
      \usepackage{amsthm}
      \usepackage{amsfonts}
      \usepackage[T1]{fontenc}{}{}
      \usepackage{latexsym}
      \usepackage{psfrag}
      \theoremstyle{plain}
      \newtheorem{thm}{Theorem}
      \newtheorem{lem}[thm]{Lemma}

      \theoremstyle{definition}
      \newtheorem{defn}{\bf Definition}
      
      \newcommand{\N}{\mathbb N}

      \newenvironment{appli}{\left( \begin{array}{ccc}}{\end{array} \right)}

      \newcommand{\ba}{\begin{appli}}
      \newcommand{\ea}{\end{appli}}
\title{A note on Hjorth's oscillation theorem.}
\begin{document}
\begin{abstract}
We reformulate an oscillation theorem proved by G. Hjorth in \cite{hjorth} in the context of continuous logic (see \cite{begnac}) and give a proof of the theorem in that setting which is similar to, but simpler than, Hjorth's original one. The point of view presented here clarifies the relation between Hjorth's theorem and first-order logic.
\end{abstract}
\maketitle

Recently, Greg Hjorth obtained a nice "oscillation theorem" for actions of Polish groups by isometries (see \cite{hjorth}). In \cite{pestov}, V. Pestov points out the importance of this result, and has this to say about its proof: ``The proof of Hjorth at this stage looks highly technical, as they say, hard. As it is being slowly digested by the mathematical community, there is no doubt that it will lead to new concepts and insights into the theory of topological groups and eventually will come to be fully understood and made into a ''soft`` proof''. This short note may be thought of as an attempt at ``digesting'' Hjorth's oscillation theorem.\\
Hjorth pointed out that his result is related to a first-order logic result (which he also proved); below we will try to understand this connection better, by proving an equivalent version of the oscillation theorem in the framework of \emph{continuous logic}. This leads to a statement mirroring the first-order one; proving the theorem in this setting also enables one to simplify the original proof a bit. Hence in a sense this note is championing the use of continuous logic (or, if not the logic, at least its language) 
to study topological groups.\\
We refer to \cite{begnac} for information about continuous logic. Below we will only deal with relational  metric structures, which we define now (we do not need to introduce here the logic of these structures).

\begin{defn}
A \emph{relational metric structure} ${\mathcal M}$ is a complete metric space $(M,d$) with $d$ bounded by $1$, along with a family $(P_i)_{i \in I}$ of uniformly continuous maps from $M^{k_i}$ to $[0,1]$ (where $k_i \in \N$ and $M^{k_i}$ is endowed, say, with the sup-metric); we always assume that the distance function $d \colon M^2 \to [0,1]$ is included in our list of predicates.\\
The structure is said to be \emph{Polish} if the underlying metric space is, that is, if $M$ is separable.\\
We say that two tuples $(a_1,\ldots,a_n)$ and $(b_1,\ldots,b_n)$ in $M^n$ have the same quantifier-free type if for all $\{j_1,\ldots,j_k\} \subseteq \{1,\ldots,n\}$ and all $i \in I$ with $k_i=k$ one has $P_i(a_{j_1},\ldots,a_{j_k})=P_i(b_{j_1},\ldots,b_{j_k})$
\end{defn}

A morphism from ${\mathcal M}$ to ${\mathcal M}$  is simply a from $M$ to $M$ that also preserves all the predicates (so in particular it is distance-preserving); it is an automorphism if it is also onto.
We endow the automorphism group $Aut({\mathcal M})$ of a relational metric structure ${\mathcal M}$ with the pointwise convergence topology, which turns it into a Polish group if $M$ is Polish (recall that a Polish group is a topological group whose topology is separable and completely metrizable). 

\begin{defn} We say that a relational metric structure ${\mathcal M}$ is \textit{appoximately ultrahomogeneous} if for any $n$-tuples $(a_1,\ldots,a_n)$ and 
$(b_1,\ldots,b_n)$ with the same quantifier-free type and any $\varepsilon >0$ there exists $g \in Aut({\mathcal M})$ such that $d(g(a_i),b_i) \le \varepsilon$ for all $i=1,\ldots,n$.
\end{defn}
Note that if ${\mathcal M}$ is separable and approximately ultrahomogeneous then any morphism is a pointwise limit of automorphisms (because morphisms preserve quantifier-free type).\\

If $(X,d)$ is a Polish metric space, its isometry group $Isom(X,d)$, endowed with the pointwise convergence topology, is a Polish group. We denote by $\delta$ the usual left-invariant distance on $Isom(X,d)$; if $G \le Isom(X,d)$ is a subgroup we denote by $\hat{G}$ the completion of $(G,\delta)$, which naturally identifies with a semigroup of isometric embeddings of $(X,d)$ into itself. Hjorth's oscillation theorem is the following:
 
\begin{thm} \emph{(Hjorth)} \label{hjorth}\\
Let $(X,d)$ be a complete separable metric space, and $G \le Isom(X,d)$ a group of cardinality bigger than one. Then there exists $x_0,x_1 \in X$ and uniformly continuous 
$$f \colon \overline{\{(\pi.x_0,\pi.x_1) \colon \pi \in G\}} \to [0,1]$$ 
such that for any $\rho \in \hat{G}$ there exist
$$(y_0,y_1), (z_0,z_1) \in \{\rho(\pi(x_0),\rho(\pi(x_1)) \colon \pi \in G\} $$
with $f(y_0,y_1)=0$ and $f(z_0,z_1)=1$.
\end{thm}

Note that it is enough to prove the preceding result when $d$ is bounded by $1$ and $G$ is closed in $Isom(X,d)$, that is  when $G$ is a Polish subgroup of $Isom(X,d)$.\\
The goal of this note is to establish the following version of Hjorth's theorem:

\begin{thm}  $ $\\
Let ${\mathcal M}$ be an approximately ultrahomogeneous Polish metric structure such that $|Aut({\mathcal M})| >1$. Then there exists a uniformly continuous $f \colon M^2 \to [0,1]$ and $(a_0,a_1) \in M^2$ such that for any morphism $\rho \colon {\mathcal M} \to {\mathcal M}$ one can find $(b_0,b_1)$ and $(c_0,c_1)$ in the image of $\rho^2$, both with the same quantifier-free type as $(a_0,a_1)$ and such that 
$f(b_0,b_1)=1$, $f(c_0,c_1)=0$.
\end{thm}

This statement mirrors the first-order result proved by Hjorth in \cite{hjorth} (corresponding to the case when $d$ only takes the values $0$ and $1$), and extends it to the context of metric structures.\\
Most (probably all) ideas in the proof below are already present in Hjorth's paper; however the proof limits the use of what he calls "messy approximations".\\

{\bf Proof.} As in the original proof, we divide the proof in subcases. In the following we let $G=Aut({\mathcal M})$. Recall that by approximate ultrahomogeneity any morphism of ${\mathcal M}$ is a pointwise limit of elements of $G$.\\

{\bf Case I.} Any $a \in M$ has a precompact orbit under $G$. Since any morphism induces an isometry of $\overline{G.a}$ into itself, and self-isometries of compact metric spaces are necessarily onto, we see that in this case any morphism is onto. Thus there is essentially nothing to prove in this case.\\

In what follows, we fix some $a$ such that $G.a$ is not precompact, and pick $\varepsilon >0$ such that $G.a$ contains infinitely many disjoint open balls of diameter $10 \varepsilon$.\\
For any $\delta>0$,  we let $ Stab_{\delta}(a)=\{g \in G \colon d(g.a,a) < \delta\}$, and
$$ acl_{\delta}(a)=\{y \in G.a \colon Stab_{\delta}(a).y \mbox{ is covered by finitely many balls of radius }\varepsilon  \}$$ 
$ $\\
{\bf Case II.} There exists $\delta >0$ such that $acl_{\delta}(a)$ is not precompact.\\
Then there is some $\tilde{\varepsilon}$ such that 
$acl_{\delta}(a)$ contains infinitely many disjoint balls of radius $\tilde{\varepsilon}$; without loss of generality we can assume $\tilde{\varepsilon}=\varepsilon$. From now on fix some countable dense $\{a_i\}_{i \in \N}$ in $G.a$.

\begin{lem}
We can find sequences $d_i,e_i$ such that $d_i,e_i \in acl_{\delta}(a_i)$ and
$$d(Stab_{\delta/2}(a_i).d_i,Stab_{\delta/2}(a_j).e_j) \ge \epsilon \mbox{ for any } i,j \in \N $$
\end{lem}

{\bf Proof of the Lemma.} Assume we have been able to define $d_i,e_i$ up to some $n$.
One can find infinitely many disjoint balls $B(z_j,10 \varepsilon)$ in $acl_{\delta}(a_{n+1})$; we need to find some $j$ such that 
$Stab_{\delta/2}(a_{n+1}).z_j$ is at distance larger than $\varepsilon$ from a set that is covered by a finite number of open balls of radius $\varepsilon$. If this is not possible, then there is an infinite $J \subset \N$ such that, for all $j \in J$, $z_j$ is mapped by some $g_j \in Stab_{\delta/2}(a_{n+1})$ at distance stricly less than $\varepsilon$ from one of these balls; so for $j,k \in J$ we get 
$$d(g_j(z_j),g_k(z_k) < 4 \varepsilon \ .$$
Fix some $j \in J$; we have $d(g_k^{-1}g_j(z_j),z_k) < 4 \varepsilon$ and from $d(z_k,z_l) >10 \varepsilon $ we obtain, for any $l\ne k \in J$:
$d(g_k^{-1}g_j(z_j), g_l^{-1}g_j(z_j))>2 \varepsilon$. Since each $g_k^{-1}g_j$ belongs to $Stab_{\delta}(a_{n+1})$, this contradicts the fact that $z_i \in acl_{\delta}(a_{n+1})$.\\
Hence one can find some suitable $z_j$, and set $d_{n+1}=z_j$; the same line of reasoning works to obtain $e_{n+1}$. This concludes the proof of the lemma . \hfill $\square$ \\
 
Now it is easy to conclude: set $D=\bigcup Stab_{\delta/2}(a_i).d_i$, $E=\bigcup Stab_{\delta/2}(a_i).e_i$. From the lemma we get $d(D,E) \ge \varepsilon$, and so it is easy to find a uniformly continuous map $f \colon M \to [0,1]$ such that $f(x)=1$ whenever $d(x,D) < \varepsilon/10$ and $f(x)=0$ whenever $d(x,E)<\varepsilon/10$.\\
Then for any morphism $\rho$ of ${\mathcal M}$ we may assume (up to multiplying $\rho$ on the right by some automorphism, which does not change the image of $\rho$) that there is some $i$ such that $d(a_i,\rho(a_i)) < \delta/2$; so in particular 
$\rho(d_i) \in \overline{D}$ while $\rho(e_i) \in \overline{E}$, and so $f(\rho(d_i))=1$ while $f(\rho(e_i))=0$.\\
(So in this case, as in case I, we obtain a function of \textit{one} variable which oscillates on the image of any morphism of ${\mathcal M}$).\\

{\bf Case III.} For any $\delta >0$ $acl_{\delta}(a)$ is precompact.\\
Following Hjorth, we pick $(z_i)$ dense in $G.a$, and find a uniformly continuous $f \colon M^2 \to [0,1]$ such that, for all $n$, $f$ equals $1$ on $B(z_n,\varepsilon/10) \times (X\setminus \bigcup_{m \le n} B(z_m,\varepsilon/2))$ while $f$ equals $0$ on $(X\setminus \bigcup_{m < n} B(z_m,\varepsilon/2)) \times B(z_n,\varepsilon/10)$.\\

\begin{lem} \label{tec2} For any $\delta >0$ there exist $(a_0,a_1)$ such that $Stab_{\delta}(a_i).a_j$ contains infinitely many disjoint open balls of radius $\varepsilon$.
\end{lem}
{\bf Proof of Lemma \ref{tec2}.} Fix $\delta > 0$. There is some $N$ such that for any $b \in G.a$ $acl_{\delta}(b)$ is covered by $N$ balls of radius $\varepsilon$; we can find $b_0,\ldots,b_N \in G.a$ such that $d(b_i,b_j)> 10 \varepsilon$ and then pick some $c \in G.a$ with 
$d(c, \bigcup_{i=0,\ldots,N} acl_{\delta}(b_i))> 5 \varepsilon $. In particular, $b \not \in acl_{\delta}(b_i)$ for all $i$; but then since 
$acl_{\delta}(c)$ is covered by $N$ balls of radius $\varepsilon$, there has to be some $i_o$ such that $b_{i_0} \not \in acl_{\delta}(c)$. \\
We just obtained $a_0$, $a_1$ both in $G.a$ and such that $a_i$ does not belong to $acl_{\delta}(a_j)$, which proves the lemma. \hfill $\square$\\

Now pick $a_0$, $a_1$ as above for $\delta=\varepsilon /20$. We claim that $(f, a_0,a_1)$ satisfy the conclusion of the theorem. Pick a morphism $\rho$ of ${\mathcal M}$ ; we can find $z_{m_0}, z_{m_1}$ such that $d(\rho(x_0),z_{m_0})< \delta$ and 
$d(\rho(x_1),z_{m_1})<\delta$. Let $k=m_0+m_1+1$. Using the lemma we can find $(\pi_i)_{i=1\ldots k} \in Stab_{\delta}(a_0)$ 
and $(\pi'_i)_{i=1\ldots k} \in Stab_{\delta}(a_1)$ such that the balls $B(\pi_i(a_1),\varepsilon)$ are disjoint, and similarly for 
$B(\pi'_i(a_0),\varepsilon)$.\\
But then we obtain that $d(\rho \circ \pi_i(a_0),z_{m_0}) < \varepsilon/10$ for all $i=1,\ldots,k$ while 
$d(\rho \circ  \pi_i(a_1), d(\rho \circ  \pi_j(a_1))> \varepsilon$ for any $i \ne j$. Hence any $z_j$, $j \le k$, can only belong to one 
ball $B(\pi_i(a_1),\varepsilon/2)$, so there is some $i_0 <k$ such that no $z_j$ belongs to $B(\pi_{i_1}(a_1),\varepsilon/2)$. Looking at the definition of $f$, we obtain $f(\rho \circ \pi_{i_0}(a_0),\rho \circ \pi_{i_0}(a_1))=1$. Similarly one finds some $j_0$ such that 
$f(\rho \circ \pi'_{j_0}(a_0),\rho \circ \pi'_{j_0}(a_1))=0$, which concludes the proof of the theorem. \hfill $\square$\\

Note that our theorem is a direct consequence of Hjorth's oscillation theorem (Theorem \ref{hjorth} above). The converse is true: let $X$ be a Polish metric space 
and $G \le Isom(X,d)$ be Polish; replacing $d$ by $d/(1+d)$ if necessary (which doesn't change either the isometry group or the uniformly continuous functions on $X$) we can assume that $d$ is bounded by $1$. Then for any $n$ consider the closed equivalence relation $\sim_n$ coming from the diagonal action of $G$ on $X^n$: $x=(x_1,\ldots,x_n) \sim_n y=(y_1,\ldots, y_n)  \Leftrightarrow x \in \overline{G.y}$ (this is an equivalence relation because $G$ acts on $X^n$ by isometries). For any $\sim_n$-class C, add a predicate 
$P_{C} \colon M^n \to [0,1]$ defined by 
$$P_{C}(x_1,\ldots,x_n) =\min\big(1,d((x_1,\ldots,x_n),C)\big) $$
We claim that the metric structure ${\mathcal M}$ obtained by adding all those predicates to $X$ is approximately ultrahomogeneous and has $G$ as its automorphism group. It is clear that any element of $G$ preserves all our predicates, and hence is an automorphism of ${\mathcal M}$; given the predicates we chose, it is then obvious that ${\mathcal M}$ is approximately ultrahomogeneous.\
To show that $G=Aut({\mathcal M})$, let $\pi$ be an automorphism of ${\mathcal M}$. Then, for any $x_1,\ldots,x_n \in X$, 
$(x_1,\ldots,x_n)$ and $(\pi(x_1),\ldots,\pi(x_n))$ have the same type, and so for all $\varepsilon >0$ there is $g \in G$ such that 
$d(g(x_i),\pi(x_i))< \varepsilon$. This is enough to show that $\pi$ is a pointwise limit of elements of $G$ and so belongs to $G$ (since $G$ is Polish it must be closed in $Isom(X,d)$). Then applying our theorem to ${\mathcal M}$ one recovers Hjorth's oscillation theorem. \\
Note that the reasoning above also shows the following result, some variants of which were already known (see for example theorem 2.4.5 in \cite{Gao})
\begin{thm}
Any Polish group is isomorphic to the automorphism group of some approximately ultrahomogeneous Polish metric structure.
\end{thm}
Actually, our technique above shows that any action by isometries of a Polish group $G$ on a Polish metric space $X$ can be seen as the action of $Aut({\mathcal X})$ on ${\mathcal X}$, where ${\mathcal X}$ is some approximately ultrahomogeneous relational Polish metric structure (with universe $(X,d)$).


\begin{thebibliography}{2}
\bibitem{begnac} I. Ben Yaacov, A. Berenstein, W. Henson and A. Usvyastov,  \textit{Model theory for metric structures}, Model Theory with Applications to Algebra and Analysis, volume 2, (Zo\'e Chatzidakis, Dugald Macpherson, Anand Pillay and Alex Wilkie, eds.), London Math Society Lecture Note Series, 350 (2008), 315-427. 
\bibitem{hjorth} G. Hjorth, \textit{An oscillation theorem for groups of isometries}, Geometric and Functional Analysis, 18(2),  Birkh\"auser (2008), 489-521.
\bibitem{Gao} S. Gao, \textit{Invariant Descriptive Set Theory}, Pure and Applied Mathematics (Boca Raton), CRC Press (2009).
\bibitem{pestov} V. Pestov, \textit{Dynamics of Infinite-dimensional Groups and Ramsey-type Phenomena}, Publica\c c \~{o}es dos Col\'oquios de Matem\'atica, IMPA, Rio de Janeiro, 2005, p. 211.
\end{thebibliography}
\end{document}